\documentclass[12pt]{article}
\font\bb=msbm10 at 12pt
\font\bbi=msbm9

\def\Z{\hbox{\bb Z}}

\def\R{\hbox{\bb R}}

\def\T{\hbox{\bb T}}

\def\Ri{\hbox{\bbi R}}

\newtheorem{Th}{Theorem}[section]
\newtheorem{pro}[Th]{Proposition}
\newtheorem{lem}[Th]{Lemma}

\newtheorem{rem}[Th]{Remark}

\def\1{1\hspace{-1.2mm}\mbox{{\normalsize I}}}

\def\Dem#1{\vskip 4mm\noindent {\bf Proof #1: }}

\title{New Proofs of  Davies-Simon's  Theorems  about 
Ultracontractivity  and Logarithmic Sobolev Inequalities
 related to Nash Type Inequalities}

\author{MAHEUX Patrick
\footnote
{Research partially supported by the European Commission (IHP  Network " Harmonic 
Analysis and Related Problems" 2002-2006, contract HPRN-CT-2001-00273-HARP)}.\\
 \thanks{D\'epartement de Math\'ematiques - Universit\'e
d'Orl\'eans - BP 6759 - 45067 ORLEANS CEDEX 2 - France. E-mail 
:pmaheux@univ-orleans.fr}}
\date{September 4, 2006}

\begin{document}
\maketitle
\begin{center}
ccsd Version
\end{center}

\begin{abstract}
We present  new proofs of two theorems of E.B. Davies and B. Simon  
 (Thm. 2.2.4  and Cor. 2.2.8 in \cite{d}) about
ultracontractivity  property  ($Ult$ for  short) of semigroups of operators and logarithmic Sobolev inequalities with 
parameter ($LSIWP$ for  short) satisfied by the generator of the semigroup. In our proof,
we    use neither the 
$L^p$ version of the $LSIWP$ 
(Theorem \ref{gammap}) nor  Stein's interpolation.
Our tool is Nash type inequality ($NTI$ for  short) as an intermediate step between 
$Ult$ and $LSIWP$.
We also  present   new results.
First, a new  formulation about the implication $LSIWP$  $\Rightarrow$ $Ult$ using a result of T.Coulhon (\cite{c}). Second,
we  show that $LSIWP$ and $NTI$ are equivalent.
We discuss   different approaches to  get Nash type 
inequalities from  an ultracontractivity property. We give  some examples of one-exponential   and double-exponential ultracontractivity and also discuss the general theory for the second case.
\end{abstract}

\\

\\
parameter, Nash type inequality,
semigroups of operators, Dirichlet form, Heat kernel, infinite Torus.
\noindent

\tableofcontents

\section{Introduction}\label{intro}
\setcounter{equation}{0}
In this paper, we give new proofs of theorems due to E.B.Davies and B.Simon (see \cite{d}
 Thm. 2.2.4  and Cor. 2.2.8 (see also  \cite{ds}). They proved  that an ultracontractivity property
 ($Ult$ for short),
\begin{equation}\label{ultdebut} 
 \vert \vert  T_t f   \vert \vert_{\infty}  \leq e^{M(t)}  \vert \vert   f  \vert \vert_1, \quad \forall t>0 \quad(Ult)_M
\end{equation}
 of a semigroup $T_t=e^{-tA}$ (under some additional assumptions) implies
 a  logarithmic Sobolev inequality with
parameter ($LSIWP$ for short),
\begin{equation}\label{lsidebut }
 \int f^2\log( f/\vert \vert   f  \vert \vert_2)\,d\mu \leq
 t(Af,f)+{\beta}(t)\vert \vert   f  \vert \vert_2^2, \quad \forall t>0 \quad (LSIWP)_{\beta}
 \end{equation}
 satisfied by the generator $A$ with $\beta=M$. They also proved some converse results with  some additional assumptions on the function $\beta$. Indeed,
it is not always true that a $(LSIWP)_{\beta}$  satisfied by a generator $A$ implies an ultracontractivity property  $(Ult)_M$ of the corresponding semigroup $T_t=e^{-tA}$. In  \cite{ds} 
(see Rmk 1 p.359), the authors give an  example of generator $A$ satisfying $(LSIWP)_{\beta}$
with $\beta(t)=ce^{1/t}$ but with no ultracontractivity. So the converse implication doesn't hold in general. But under some conditions on the function $\beta$, it can be proved that $(LSIWP)_{\beta}$ implies  $(Ult)_{\tilde{M}}$
 with some function ${\tilde{M}}$ (which may differ of the function $M$ in (\ref{ultdebut})). An interesting situation is when 
  ${\tilde{M}}(t)=c_1M(c_2t)+c_3$.  In that case, the two statements $(LSIWP)_{\beta}$ and $(Ult)_M$  are equivalent in the sense that $M$ and ${\tilde{M}}$ behave in the same way.
 For example,  $e^{M(t)}= Ce^{-\lambda t} t^{-d}e^{c/t^{\gamma}}$
i.e  $M(t)=k_1-\lambda t-d\ln t +c/t^{\gamma}$ with $ k_1=\ln C,\lambda,d,c,\gamma>0$. But we are unable to prove this relation between 
 ${\tilde{M}}(t)$ and $M(t)$ in the general situation.
 For instance, if $\beta(t)=e^{c/t^{\gamma}}$ with $0<\gamma<1$, we are only able to prove that 
 ${\tilde{M}}(t)=e^{c/t^{{\gamma}'}}$ with ${\gamma}'=\frac{\gamma}{1-\gamma}$ as far as the author knows.
 It could be conjectured that  ${\gamma}'=\frac{\gamma}{1-\gamma}$ is optimal. In particular, the singularity of the behavior of  ${\gamma}'$  clearly appears when $\gamma$ goes to 1. 
 \\
 
 Let us recall briefly the interest of LSIWP and   ultracontractivity property.
  The ultracontractivity  property is equivalent for (symmetric) semigroup to the following on-diagonal estimate of the heat kernel:
  \begin{equation}\label{heatdebut }
 \sup_x h_t(x,x)=\sup_{x,y}h_t(x,y)  \leq e^{M(t)}, \quad \forall t>0 
\end{equation}
 with
   \begin{equation}
  T_t f (x)=   \int  h_t(x,y)\,d{\mu}(y)
\end{equation}
($h_t(x,y)$   is the so-called  heat kernel). 
\\

So, when we are able to prove a LSWIP and we have at hand a theorem saying that LSIWP implies
an ultracontractivity property with an explicite bound, we  immediately deduce  the same bound on the heat kernel.
 Note  that,  if we replace LSIWP by a Nash type inequality, the same remark holds true (see \cite{c}). 
 Recall that a Nash type inequality (NTI for short)
  is the control by a function  $\Theta$ of the $L^2$-norm by the quadratic form associated to the generator of the semigroup when the $L^1$-norm is bounded. More precisely,  
     \begin{equation}
 \Theta\left(   \vert\vert  f\vert\vert_2^2\right)
  \leq
  (Af,f) ,\quad \forall f\in {\mathcal D}(A),\;   \vert\vert 
f\vert\vert_1\leq 1 .
  \end{equation}
  
  Since LSIWP and Nash type inequality have the main goal (To prove ultracontractivity), we may ask for relations  between these two inequalities.  In that paper, we show that  LSIWP and Nash type inequality are equivalent. Moreover, we use Nash type inequality  to give a new proof of the statement $(LSIWP)$ 
  $\Rightarrow$  $(Ult)$. We   use  again Nash type inequality  for the converse implication. This  sheds new lights on relationships between these three inequalities.
\\

We now describe the contents  of this paper.
\\

 In Section 2, we describe well-known results about   the relationship
between   ultracontractivity property and logarithmique Sobolev inequalities with parameter following
\cite{ds}. We also recall  part of the proofs  for the convenience of the reader and also to be compared  with the new methods developed  in that paper.
\\

In Section 3, we prove the implication 
$(LSIWP)_{\beta}$ $\Rightarrow$ $(Ult)_M$
under the usual 
additional assumptions on  the function $\beta$.  We use a new approach for the proof:
we introduce an intermediate step with Nash type inequality.
 In fact, we give two different results. The first result is  a new proof of a corollary   of a general  result of Davies and Simon (see \cite{ds}). The second proof gives another ultracontractive bound using a result of T.Coulhon \cite{c}.
This last result  has more general applications.  
\\

In Section 4, we study the converse:
$(Ult)_M$  $\Rightarrow$ $(LSIWP)_{\beta}$. The new proof has two steps.
We deduce a  Nash type inequality from the
ultracontractivity  property by using again a result of (\cite{c}). This  Nash type inequality 
is equivalent to a Nash type inequality with parameter already close to $(LSIWP)_{\beta}$.  The second step consists in  applying  a method of truncation for Dirichlet forms developed in \cite{bcls}
to obtain $(LSIWP)_{\beta}$.
\\

In Section 5, we prove the equivalence between LSIWP and Nash type inequality using ideas developed in the preceding sections.
\\

In Section 6, we discuss different well-known approaches  to prove Nash type inequalities. Such discussion arises naturally  since  Nash type inequality is the main tool  of our  proofs.
\\

In Section 7, we focus on the polynomial ultracontractivity  property  i.e $e^{M(t)}=ct^{-\nu}, \, \nu>0$ in (\ref{ultdebut}) which is equivalent to the usual $L^2$-Sobolev inequality. We show how this  ultracontractivity  property can be expressed in terms of different  functional inequalities. In particular, we recall the weak-Sobolev inequality introduced earlier by D. Bakry (see \cite{b}). This section essentially collects these information.
\\

Section 8. To show how this general theory can be applied outside of the polynomial ultracontractivity  property setting, we mention two families of examples of heat kernels on the infinite dimensional torus $\T^{\infty}$ coming from \cite{b2}.  The generator of the semigroup is an infinite dimensional Laplacian with constant coefficients. The sequence of these  coefficients has to go to infinity. Depending on  the speed of this sequence, the corresponding heat kernel bound  has a different behavior. The first family of examples is the one-exponential ultracontractivity  behavior i.e
$$
e^{M(t)}=c_1e^{c_2/t^{\gamma}},\quad (\gamma >0).
$$
in (\ref{ultdebut}).
In this situation,  $(Ult)_M$  and the corresponding $(LSWIP)_{\beta}$  are equivalent (for  the abstract theory) in the sense that 
$M(t)=c/t^{\gamma}$ and ${\beta}(t)=c/t^{\gamma}$ with possibly different constants $c$.

The second familly of examples is the double-exponential ultracontractivity  behavior i.e
$$
e^{M(t)}=c_1e^{e^{c_2/t^{\gamma}} },\quad  (\gamma >0).
$$
in (\ref{ultdebut}).
In this case, the general theory doesn't give equivalence between  $(Ult)_M$  and the corresponding $(LSWIP)_{\beta}$.
More examples of heat kernel behaviors (i.e ultracontractivity  ) can be found in \cite{b2}. 
\\

In Section 9, we discuss in the general frame work   the relationship between double-exponential ultracontractivity  property
and the corresponding  $(LSWIP)_{\beta}$ (or the corresponding Nash type inequality).

\section{Relations between ultracontractivity and LSIWP}\label{relultlsi}
\setcounter{equation}{0}
Let 
$(e^{-At})_{t\geq 0}= (T_t)_{t\geq 0}$ be a symmetric Markov 
semigroup on $ L^2(X,d{\mu})$ 
with generator $A$ defined on   a ${\sigma}$-finite 
measure space $(X,d{\mu})$.
We say that $(T_t)_{t\geq 0}$ is ultracontractive if 
for any $t>0$, there exists a finite positive number $a(t)$ such 
that, for all $f\in L^1$ :
\begin{equation}\label{ult1}
\|T_tf\|_{\infty}  \leq a(t) \|f\|_1.
\end{equation}

An equivalent formulation (by interpolation) of ultracontractivity is 
as follows:
For any $t>0$, there exists a finite positive number  $c(t)$ such 
that, for all $f\in L^2$,
\begin{equation}\label{ult2}
\|T_tf\|_{\infty} \leq c(t) \|f\|_2
\end{equation}
  Also by duality, the inequality (\ref{ult2}) is equivalent to
\begin{equation}\label{ult3}
\|T_tf\|_{2} \leq c(t) \|f\|_1
\end{equation}
It is known that, under the assumptions on the semigroup 
$(T_t)_{t\geq 0}$, (\ref{ult2}) implies (\ref{ult1})
with $a(t)\leq c^2(t/2)$
and 
(\ref{ult1}) implies (\ref{ult2})  with $c(t) \leq \sqrt{a(t)}$.
\\

We say that the generator $A$ satisfies  LSIWP  (logarithmic Sobolev inequality 
with parameter) if 
 there exist a monotonically decreasing continuous function 
${\beta}: (0,+\infty)\rightarrow (0,+\infty)$ such that
\begin{equation}\label{lsiwp}
\int f^2\log f\, d{\mu} \leq
\epsilon Q(f) +{\beta}(\epsilon) \|f\|^2_2 + \|f\|^2_2\log \|f\|_2
\end{equation}
for all $\epsilon >0$ and $0\leq f\in \mbox{Quad}(A)\cap L^1\cap 
L^{\infty}$ where
$\mbox{Quad}(A)$ is the domain of $\sqrt{A}$ in $L^2$ and 
$Q(f)=(\sqrt{A}f,\sqrt{A}f)$.
\\

This inequality is modeled on the celebrated Gross inequality \cite{g}.
\\

In \cite{ds},\cite{d}, the authors show that LSIWP implies 
ultracontractivity property  under an integrablity condition on 
$\beta$. This condition can be enlarged and be stated as follows:

\begin{Th}\label{coroda}(Cor. 2.2.8 \cite{d} ).
Let ${\beta}(\epsilon)$ be a monotonically decreasing continuous 
function of $\epsilon$
such that
\begin{equation}\label{vareps}
\int f^2\log f \, d{\mu}\leq
\epsilon Q(f) +{\beta}(\epsilon)\, \|f\|^2_2 + \|f\|^2_2\log \|f\|_2
\end{equation}
for all $\epsilon >0$ and $0\leq f\in \mbox{Quad}(A)\cap L^1\cap 
L^{\infty}$. Suppose that
for one ${\eta}>-1$,
\begin{equation}\label{integral}
M_{\eta}(t)=({\eta}+1)t^{-({\eta}+1)})\int_0^t 
{s}^{\eta}{\beta}\left(\frac{s}{\eta+1}\right)
\,ds
 \end{equation}
is finite for  all $t>0$. Then $e^{-At}$ is ultracontractive 
and 
\begin{equation}\label{majo}
\| e^{-At} \|_{\infty,2}\leq e^{M_{\eta}(t)}
\end{equation}
for all $0<t<\infty$.
\end{Th}

Before recalling the proof of   Davies and Simon, we make some 
comments.
\\

Corollary  
2.2.8 of \cite{d}  is Theorem \ref{coroda} with $\eta=0$. 
In  the literature,  Corollary  2.2.8  of \cite{d} is used, for instance,  to deal 
with   $\beta(t)= ct^{-\alpha}$  for $0<\alpha <1$ and we had to go back to Theorem  2.2.7 to deal with the case $\alpha\geq 1$ (see \cite{fl}  for such an instance of application).
 Theorem  \ref{coroda} unify these two cases in one case
 just by an 
appropriate choice of $\eta$.
Indeed it is easy to obtain the bound of ultracontractivity in the theorem above with the parameter $\eta$ by 
the  same argument used to treat the example 2.3.4 p.72 of  \cite{d}.
The proof  will be 
recalled below. Note that, in general for our applications, $\beta$ is non-increasing so that for any $\eta>-1$, we have  ${\beta}(\frac{t}{\eta+1})\leq M_{\eta}(t), t>0$.
 \\

So the interest of such result relies on the  fact that we can choose 
the parameter  $\eta $. Indeed for some parameter $\eta$ the integral 
(\ref{integral}) may not converge at the origin but   it may
converge for some other parameters $\eta$. 
For instance when $\beta(t)=c/t^{\alpha}$ ($\alpha>0$), we obtain
$M_{\eta}(t)=c'/t^{\alpha}$ with the same index $\alpha$ but we have to choose
$\eta>\alpha-1$. 
The weight $s^{\eta}$   is used to remove the singularity of the integral at 
the origine. So for  this example of class of functions, with  an appropriate choice of  $\eta$, the 
integral (\ref{integral}) converges and we recover the 
function $\beta$ (up to a multiplicative constant).  

It may also happened that, for some function $\beta$, the integral doesn't converge for any choice of $\eta$. For instance,  $\beta(t)=\exp(c/t^{\alpha}), \alpha>0$.
\\

  The  aim  of this paper is to give a different proof of 
this result (see Section \ref{secproof}).
\\

We now recall the main steps of the proof of Theorem  \ref{coroda}
for the case $\eta=0$ and give the  proof of the general case 
$\eta >-1$ (which can be deduce from Example 2.3.4 p.72 of \cite{d}).

The first step   is the following lemma. 
This 
lemma says that if  an $L^2$-version of  LSIWP is satisfied  then an  $L^p$-version is also satisfied
 for any $p\in (2,+\infty)$.

\begin{lem}\label{deducgam}(Lemma 2.2.6 \cite{d})
Assume that the LSIWP (\ref{lsiwp}) is satisfied. Then
\begin{equation}\label{lp}
\int g^p\log g \, d{\mu}\leq
{\epsilon}(Ag,g^{p-1})+ 2{\beta}({\epsilon})p^{-1}\, \|g\|^p_p +
\|g\|^p_p\log \|g\|_p
\end{equation}
for all $2<p<\infty$, all ${\epsilon}>0$  and all $g\in {\cal 
D}_+=\bigcup_{t>0}e^{-At}(L^1\cap L^{\infty})_+$.
\end{lem}
For the next step, the parameter $\epsilon$ can be
chosen as a function of $p$  in the $L^p$-inequality (\ref{lp}).  Then we can deduce the 
ultracontractivity property from this family of $L^p$-inequalities. 

\begin{Th}\label{gammap}(Thm 2.2.7 \cite{d})
Let ${\epsilon}(p)>0$ and ${\Gamma}(p)$ be two continuous functions 
defined for $2<p<\infty$
such that
\begin{equation}
\int f^p\log f\, d{\mu} \leq
{\epsilon}(p) <Af,f^{p-1}>+ \,{\Gamma}(p)\,\|f\|^p_p +
\|f\|^p_p\log \|f\|_p
\end{equation}
for all $2<p<\infty$ and all $f\in {\cal 
D}_+=\bigcup_{t>0}e^{-At}(L^1\cap L^{\infty})_+$.
If  
\begin{equation}
t=\int_2^{\infty}p^{-1}{\epsilon}(p)\, dp, \quad 
M=\int_2^{\infty}p^{-1}{\Gamma}(p)\, dp
\end{equation}
are both finite then $e^{-At}$ maps $L^2$ into $L^{\infty}$ and 
\begin{equation}
\| e^{-At} \|_{\infty,2}\leq e^M
\end{equation}
\end{Th}

{\bf Proof of Theorem \ref{coroda}:} Let $\eta>-1$ and set 
${\nu}^{-1}={\eta}+1>0$.
We apply   Lemma \ref{deducgam}  and  Theorem \ref{gammap}  with 
$${\epsilon}(p)=t{\nu}2^{\nu}p^{-{\nu}},\qquad 
{\Gamma}(p)=2{\beta}({\epsilon}(p))p^{-1},$$
then 
$M(t)=\int_2^{+\infty} 2{\beta}({\epsilon}(p))p^{-2}\,dp=M_{\eta}(t)$ 
as defined above.
 This completes the proof.
\\


Now we consider  the converse statement.
We recall the following result due to Davies and Simon. In this statement, we note that there is no restriction as the integrability condition of Theorem \ref{coroda}. Thus an ultracontractivity property always implies LSIWP.

\begin{Th}\label{gross}(Thm 2.2.3 \cite{d})
Assume that $e^{-At}$is  ultracontractive i.e.
\begin{equation}\label{ultraM}
\| e^{-At} \|_{\infty,2}\leq e^{M(t)}
\end{equation}
for all $t>0$, where $M(t)$ is  a monotonically decreasing continuous 
function of $t$.

Then $0\leq f\in \mbox{Quad}(A)\cap L^1\cap L^{\infty}$ implies 
$f^2\log f\in L^1$,
and the logarithmic Sobolev inequality 
\begin{equation}\label{sleps}
\int f^2\log f \, d{\mu}\leq
\epsilon Q(f) +M(\epsilon) \|f\|^2_2 + \|f\|^2_2\log \|f\|_2
\end{equation}
for all $\epsilon>0$.
\end{Th}
\Dem{}
We just recall the main arguments.  We consider
$Q_z=e^{-tzA}$ for $0\leq \Re z \leq 1$ for a fixed $t$.
For any  $y\in \R$, we have 
$$
\vert\vert Q_{iy}f\vert\vert_{2}\leq \vert\vert  f \vert\vert_2
$$
and 
$$
\vert\vert Q_{1+iy}\vert\vert_{\infty}\leq e^{M(t)} \vert\vert  f \vert\vert_2.
$$
By Stein's complex interpolation Theorem, for any $0<s<t$,
with $\theta=s/t$, we have
$$
\vert\vert Q_{\theta}\vert\vert_{p(s)}\leq e^{\theta M(t)}\vert\vert  f \vert\vert_2
=e^{sM(t)/t}\vert\vert  f \vert\vert_2
$$
with $p(s)=2t/(t-s)$. Note that at this stage, the dependance of the bound in $s$ is very simple and  have the value $1$ at $s=0$. 
\\

The second idea is to obtain the expression under the integral with $\ln f$ by deriving  at $s=0$  the $L^{p(s)}$-norm of $f_s=T_sf$ (here we skip  the details).
Let  $\phi(s)=\vert\vert  f_s \vert\vert_{p(s)}^{p(s)}$,
$$
\phi'(s)
=
p(s)<-Af_s,f_s^{p(s)-1}>+\,p'(s)\int f_s^{p(s)}\ln f_s \, d{\mu}.
$$
Let 
$\psi(s)=e^{sM(t)/t}$ and assume $\vert\vert  f \vert\vert_2=1$ thus $\phi(0)=\psi(0)=1$!
Consequently,
$$
\phi'(0)\leq \psi'(0)
$$
that is
$$
-2<Af,f>+\frac{2}{t}\int f^2\ln f\, d{\mu} \leq 2M(t)/t.
$$
The proof is completed.
\\

 Our new proof in Section \ref{nproflsi} avoid the interpolation argument .
 \\
 
The main applications we have in mind for  Theorem \ref{coroda} is 
with
$${\beta}({\epsilon})=\ln(a(\epsilon))=
\ln c_1 -{\lambda}{\epsilon}-d\ln {\epsilon}+c_2/{\epsilon}^{\gamma}$$
with $c_1,c_2>0$ and ${\lambda}, d, 
{\gamma}\geq 0$ (i.e.
$a(t)=c_1e^{-{\lambda}t}
{t}^{-d}\exp(c_2/{t}^{\gamma})$). For a suitable choice of
${\eta}$ in Theorem \ref{coroda}, we obtain for $M_{\eta}$ a function of the same type as 
${\beta}$.  More precisely
$M_{\eta}(t)=\ln c_1'-\lambda't-d\ln t+c_2'/{t}^{\gamma}$. So 
$e^{M_{\eta}(t)}=c_1'\exp(-\lambda't)t^{-d}\exp(c_2'/{t}^{\gamma})$ is of the same type as $a(t)$ above (that is up to constants $c_1$, $\lambda$, $c_2$). Note that the exponents $d$ and $\gamma$ are preserved in this transformation.
 For this class of function,  Theorem \ref{coroda}
and Theorem \ref{gross} are converse of each other.
\\

Of course, other classes  of functions $\beta$ can be considered but we are not always able to pass from 
LSIWP  to  ultracontractivity. Indeed, it is worth noting that there exists also a semigroup which is not ultracontractive 
but satisfies   LSIWP (\ref{vareps})
with  ${\beta}({\epsilon})=c_1\exp(c_2/{\epsilon})$ (see p.359 and also Section 4 p.355 and section 5 p.357
 in \cite{ds}).
At the end of this paper, we give an alternative proof for an explicit bound for 
ultracontractivity when  ${\beta}(\epsilon)=c_1\exp(c_2/{\epsilon}^{\gamma})$, ${\gamma}\in ]0,1[$.
\\

To finish this section, we mention the following result.
If we only suppose $0\leq f\in Quad(A)$ then it is not obvious that
$f^2\log f \in L^1$ (at least when $\mu$ is not finite) and  we have a slight variation of Theorem
\ref{gross}.

\begin{Th}\label{logplus}(Thm 2.2.4 \cite{d}).
Let $e^{-At}$ be ultracontractive with 
\begin{equation}
\| e^{-At} \|_{\infty,2}\leq e^{M(t)}
\end{equation}
for all $t>0$, where $M(t)$ is  a monotonically decreasing continuous 
function of $t$.
Then $0\leq f\in \mbox{Quad}(A)$ implies 
\begin{equation}
\int f^2\log_+ f \, d{\mu}\leq
\epsilon Q(f) +{\beta}(\epsilon) \|f\|^2_2 + \|f\|^2_2\log \|f\|_2
\end{equation}
for all ${\epsilon}>0$ where ${\beta}(\epsilon)=M({\epsilon}/4)+2$.
\end{Th}

In the next section, we give  a new
proof  of the  theorem \ref{coroda} in the form of Theorem 
\ref{logplus}. 

\section{LSWIP  implies ultracontractivity }\label{secproof}
\setcounter{equation}{0}
\subsection{Davies-Simon result}
Assuming LSIWP is satisfied by  the generator, we give a new proof of  ultracontractivity property of the associated semigroup (under some integrability conditions).  We do not use $L^p$ version of LSWIP
as in \cite{ds}.
We only use of the $L^2$ inequalities. There are three steps in our proof. First step:  from LSIWP we  deduce a (relaxed) Nash type inequality for the generator using a convexity argument (Lemma \ref{convex}). Second step: we derive a differential inequality satisfied by the associated semigroup.
Third step: we prove a universal bound on all solutions of this differential inequality (Lemma \ref{lemphi})
and, as a consequence, we deduce the 
 ultracontractivity property .

The first lemma depends on a convexity argument  (Jensen inequality). This lemma will also   be   used in section  \ref{nproflsi} and \ref{nti}.  This lemma  comes from \cite{bima}.
\begin{lem}\label{convex}
If $f\in L^1\cap L^{\infty}$ with $f\geq 0$ and $\| f\|_1=1$ 
then
\begin{equation}\label{conv}
\|f\|^2_2\log \|f\|_2\leq 
\int f^2\log (f/\|f\|_2)\, d{\mu}
\end{equation}
and, more generally, 
if $f\in L^1 \cap L^{\infty}$ with $f\geq 0$, 
\begin{equation}\label{conv2}
\|f\|^2_2\log \|f\|_2\leq 
\int f^2\log f\,d{\mu}-\|f\|^2_2\log \|f\|_2
+\|f\|^2_2\log \|f\|_1
\end{equation}
In particular,
if  $\| f\|_1\leq1$
then (\ref{conv}) holds true.
\end{lem}

\Dem {}
If $f\in L^1\cap L^2$ with $f\geq 0$ and $\| f\|_1=1$ 
then $d{\nu}=fd{\mu}$ is a probability measure. For every convex 
function
 ${\Psi}:\R^+Ê\longrightarrow \R^+$,
the Jensen inequality yields 
\begin{equation}
{\Psi}\left( \int f\,d{\nu}\right)\leq \int {\Psi}(f)\,d{\nu}
\end{equation}
We apply this to the convex function ${\Psi}(x)=x\log x$.
 Therefore, $ \int f\,d{\nu}=\| f \|_2^2$
and
\begin{equation}
2\|f\|^2_2\log \|f\|_2 \leq \int f^2\log f \, d{\mu} 
\end{equation}
We conclude (\ref{conv}) with $\| f\|_1=1$. To obtain the general 
case (\ref{conv2}), we put $f \|f\|^{-1}_1$ instead
 of $f$ 
in  (\ref{conv}) then the
inequality follows.
When $\|f\|_1\leq 1$, (\ref{conv2}) implies (\ref{conv}). The 
lemma is proved.
\begin{rem}
We can prove an analogue to the inequality (\ref{conv}) with $\log_+$ 
instead of $\log$. 
\end{rem}

The  following lemma gives  a universal bound for all the solutions of the differential inequality
with parameter (\ref{ineqphi}).
We shall note that this  bound doesn't depend on the initial condition.
We also discuss the optimality 
of the result.
\\

We need to introduce some notations and definitions.
For any ${\eta},\lambda\in \R$  and for any  continuous real-valued  
function $b$ defined 
on $(0,+\infty)$, we define  the following function:
\begin{equation}\label{hform}
H_{{\eta},\lambda, b}(t)=\frac{2\lambda}{t^{{\eta}+1}}\int_0^t 
s^{\eta}b(s/{\lambda})\, ds
\end{equation}
(assuming  this integral converges).
This function plays the role of the function $M_{\eta}$ of Section \ref{intro}.
We shall denote by $H_{{\eta},b}$ the function $H_{{\eta},(\eta+1)/2, 
b}$.

\begin{lem}\label{lemphi}
Let $t_0\in (0,+\infty]$.
Assume that ${\Phi}\in
C^1((0,+\infty),\R)$  is a   function  satisfying  the 
following differential inequality
\begin{equation}\label{ineqphi}
{\Phi}(s)\leq
(-t/2){\Phi}^{\prime}(s) +b(t) 
\end{equation}
for all $s>0$ and all $0<t<t_0$,
with $b$ a continuous real-valued  function defined on $(0, t_0)$.

\begin{enumerate}
\item
For any  ${\eta}>-1$, let $\lambda=\frac{\eta+1}{2}$.
Assume that $H_{{\eta},b}(t)$ converges for all $t\in (0,\inf (t_0,s_0/\lambda))$.
Then  for all  $t\in (0,\inf (t_0,s_0/\lambda))$:
\begin{equation}\label{ineqho}
{\Phi}(t)\leq
H_{{\eta}, b}(t).
\end{equation}

\item
Assume that $\Phi$  is a non-negative function.
For any  ${\eta}>-1$ and  $\lambda\geq \frac{\eta +1}{2}$, we 
have
 for all  $0<t<\inf (t_0,s_0/\lambda)$:
\begin{equation}\label{ineqh}
{\Phi}(t)\leq
H_{{\eta},\lambda, b}(t)
\end{equation}
Moreover, if  $b$ is non-increasing and non-negative, the function $\lambda \rightarrow H_{{\eta},\lambda, b}(t)$ is 
increasing for fixed $t$ and 
fixed $\eta$.
Hence, in that case, inequality (\ref{ineqh}) with  $\lambda=\frac{\eta+1}{2}$ implies 
the others.

\item
 For ${\eta}>-1$ and $\lambda= \frac{\eta +1}{2}$.  Then $H_{{\eta},\lambda, b}(s)$
satisfies (\ref{ineqphi}) for all $s\in (0,\inf (t_0,s_0/\lambda))$ with $t=s/\lambda$.
In fact  in that case,  the inequality (\ref{ineqphi}) is an equality.
\end{enumerate}
\end{lem}

\Dem {of Lemma \ref{lemphi}}
Let ${\eta}>-1$ and $\lambda>0$. For  $s>0$ 
choose $t=s/{\lambda}$ in (\ref{ineqphi}). We multiply 
(\ref{ineqphi}) by $s^{\eta}$
and integrate over the interval $(0,t]$. Then

\begin{equation}
\int_0^t s^{\eta}{\Phi}(s)\, ds\leq
(-1/2{\lambda})\int_0^t s^{{\eta}+1}{\Phi}^{\prime}(s)\, ds +
\int_0^t s^{\eta}b(s/{\lambda})\, ds.
\end{equation}
The second integral is integrated by parts, so we get
$$
\frac{1}{2\lambda}t^{\eta+1}\Phi(t)
\leq
\left[ \frac{\eta+1}{2\lambda}-1\right]
\int_0^t s^{\eta}{\Phi}(s)\, ds\
+
\int_0^t s^{\eta}b(s/{\lambda})\, ds.
$$
Let  $\lambda=(\eta+1)/2$ if $\Phi$ is real-valued and let
$\lambda\geq (\eta+1)/2$ if $\Phi\geq 0$. The 
second term above is negative. Then
$$
\frac{1}{2\lambda}t^{\eta+1}\Phi(t)
\leq
\int_0^t s^{\eta}b(s/{\lambda})\, ds.
$$
This   proves 1 and 2 of the lemma. The statement 3 is easy to check by a direct computation.
\\

Now, we can restate Theorem \ref{coroda} with a slight modification
 of the expression of the bound (\ref{majo}).
The proof depends upon  Lemmas \ref{lemphi}
 and \ref{convex} and it is short. This proof is also
simpler than the original proof of Theorem \ref{coroda} because it doesn't 
use Theorem
\ref{gammap}. But as already mentioned, an important disadvantage of Theorem 
\ref{coroda}  or Corollary 2.2.8 of \cite{d} is that it doesn't enable
us to treat the case 
${\beta}(\epsilon)=e^{\frac{c_1}{{\epsilon}^{\alpha}}}, {\alpha}>0$. By a
modification of   Lemma \ref{lemphi} we shall   provide in  Section  \ref{nop} an explicit  bound  of ultracontractivity  property under the assumption ${\beta}(\epsilon)=e^{\frac{c_1}{{\epsilon}^{\alpha}}}$ 
with $0<\alpha <1$ .  

\begin{Th}\label{our}
Suppose that the following logarithmic Sobolev inequality is valid 
for all $\epsilon>0$
and all $0\leq f\in \mbox{Quad}(A)\cap L^1\cap L^{\infty}$,
\begin{equation}\label{soblog}
\int f^2\log f\,d{\mu} \leq
\epsilon Q(f) +{\beta}(\epsilon) \|f\|^2_2 + \|f\|^2_2\log \|f\|_2
\end{equation}
with ${\beta}$ a continuous function.
Then for all ${\eta}>-1$ 
\begin{equation}\label{majobis}
\| e^{-At} \|_{\infty,2}\leq e^{M_{\eta}(t)}
\end{equation}
for all $0<t<\infty$ where $M_{\eta}(t)$ is defined in (\ref{integral}).
\end{Th}

\Dem{}
Assume that $0\leq f\in 
\mbox{Quad}(A)\cap  L^1 \cap  L^{\infty}$. 
We set  $f_s= e^{-As}f=T_sf$ for $s>0$. Suppose $\| f\|_1= 1$
then, by contraction property on $L^1$ of 
the semigroup,  $\|f_s\|_1 \leq \|f\|_1= 1$.
We check  that $0\leq f_s\in \mbox{Quad}(A)\cap L^1\cap L^{\infty}$.
We put $f_s$ in (\ref{soblog}),

\begin{equation}
\int f_s^2\log f_s \,d{\mu} -
\| f_s\|^2_2\log \| f_s\|_2\leq
tQ(f_s) +{\beta}(t) \|f_s\|^2_2 
\end{equation}
for all $s,t>0$.
We apply Lemma \ref{convex} with $f_s$,  we deduce
\begin{equation}
\| f_s\|^2_2\log \| f_s\|_2\leq
tQ(f_s) +{\beta}(t) \|f_s\|^2_2 
\end{equation}

Let  ${\Psi}(s)= \|f_s\|^2_2$. Then ${\Psi}^{\prime}(s)=-2Q(f_s)$. 
Therefore,
 ${\Psi}(s)$ satisfies 

\begin{equation}
\frac{1}{2}{\Psi}(s)\log {\Psi}(s)\leq
(-t/2){\Psi}^{\prime}(s) +{\beta}(t) {\Psi}(s)
\end{equation}

Let ${\Phi}(s)=\log {\Psi}(s)$ and changing $t$ by $t/2$, then
\begin{equation}
{\Phi}(s)\leq
(-t/2){\Phi}^{\prime}(s) +b(t) 
\end{equation}
for all $s,t>0$ with $b(t)=2{\beta}(t/2)$.
We apply  Lemma \ref{lemphi} with $t_0=+\infty$, for all ${\eta}>-1$ 
\begin{equation}
{\Phi}(t)=\log \| T_tf\|^2_2\leq
H_{{\eta},b}(t)=2M_{\eta}(t)
\end{equation}
Hence, 
\begin{equation}
\| T_tf\|_2\leq \exp\left( M_{\eta}(t)\right)\|f \|_1
\end{equation}
By duality,
\begin{equation}
\| T_tf\|_{\infty}\leq \exp\left( M_{\eta}(t)\right) \|f\|_2.
\end{equation}
The proof of the theorem is completed.
\\ 

We can localize this result in the sense that, if we assume (\ref{soblog}) holds true for all $0<{\varepsilon}<{\varepsilon}_0$  then  (\ref{majobis}) holds true for $t\in(0,2{\varepsilon}_0 )$.

 \subsection{Another ultracontractive bound}
 Now, we prove another ultracontractive bound for the semigroup under the assumption that the generator satisfies  LSIWP.  We use  a result due to T. Coulhon (see Prop. II.1 of \cite{c}) and the lemma (\ref{convex}) of this article.
  
 \begin{Th}\label{anotherbound}
 Let $A$ be a generator of a submarkovian semigroup and 
$\beta$ be a function such that:
\begin{equation}\label{lsiwpcoul}
 \int f^2\log( f/\vert \vert   f  \vert \vert_2)\,d\mu \leq
 t(Af,f)+\beta(t)\vert \vert   f  \vert \vert_2^2,  t>0.
 \end{equation}
We define $B(y)=\sup_{t>0}(ty/2-t\beta(1/t)), y\in \R$ and 
$M(t)$ the inverse function of $q(s)=\int_s^{+\infty}\frac{dy}{B(y)}, s\in \R$
(We assume that 
$\int^{\infty} \frac{dy}{B(y)}<\infty$ ). Then
\begin{equation}
\| T_tf\|_{\infty}\leq e^{ M(t)} \|f \|_1
\end{equation}
for any $t>0$.
 \end{Th}

 \Dem{}
 We assume that (\ref{lsiwpcoul}) is satisfied.  The first setp is to obtain a Nash type inequality. For that purpose, we apply lemma (\ref{convex}):
   \begin{equation}\label{lsiwpcoul2}
 \vert\vert  f\vert\vert_2^2  \; \ln  \vert\vert  f\vert\vert_2
  \leq
  t(Af,f)+ {\beta}(t) \vert\vert  f\vert\vert_2^2  ,\; \forall t>0, \quad \forall f\in {\mathcal D}(A),\;   \vert\vert 
f\vert\vert_1\leq 1 .
  \end{equation}
 Hence
    \begin{equation}
 \vert\vert  f\vert\vert_2^2  \left[\frac{1}{2t}\; \ln  \vert\vert  f\vert\vert_2^2-\frac{1}{t}{\beta}(t) \vert\vert  f\vert\vert_2^2 \right]^2
  \leq
  (Af,f).
  \end{equation}
 By optimisation over $t>0$ and by definition of $B$,
    \begin{equation}\label{lsiwpcoul4}
 \vert\vert  f\vert\vert_2^2  \;B\left(\ln  \vert\vert  f\vert\vert_2^2\right)
  \leq
  (Af,f) ,\quad \forall f\in {\mathcal D}(A),\;   \vert\vert 
f\vert\vert_1\leq 1 .
  \end{equation}
  Since $B$ is convex, it follows that $B$  is  continuous. 
  Let's denote by $\Theta(x)=xB(\ln x), x>0$ (here we use notations of \cite{c}).
For the second step, we apply Prop.II.1 of  \cite{c} which says that
 $$
\| T_tf\|_2\leq m(t) \|f \|_1
$$
with $m(t)$ the inverse function of $p(t)=\int_t^{+\infty} \frac{dx}{\Theta(x)}, t>0$.
By a change of variable, we get
$$
p(t)=\int_{\ln t}^{+\infty} \frac{dy}{B(y)}, t>0.
$$
 Setting 
$q(s)=\int_{s}^{+\infty} \frac{dy}{B(y)}, s\in \R$, we obtain
$p(t)=q(\ln t),t>0$. Thus $m(t)=exp(q^{-1}(t))=e^{M(t)}$ with $M(t)$ defined as in this theorem.
This completes the proof.

   
\section{Ultracontractivity implies LSIWP}\label{nproflsi}
\setcounter{equation}{0}
In \cite{d} and \cite{ds}, the authors show how we can deduce  LSIWP from 
ultracontractivity using the complex interpolation of Stein.
In this section, we give another way to deduce logarithmic Sobolev 
inequality with parameter from ultracontractivity.  
Nash type inequality  is involved in the proof.
 \\

 We first recall a result due to T.Coulhon which is one step to prove 
LSIWP from ultracontractivity property.  We give a slightly different presentation of  the statement of this theorem and we recall the proof for the convenience of the reader.
 
 \begin{Th}\cite{c}\label{coulhon}
  Let $(T_t)$ be a  symmetric semigroup on $L^2$.
  Suppose that 
   \begin{equation}\label{ultracoul}
  \vert\vert T_t f\vert\vert_2^2\leq m(t)
   \vert\vert f\vert\vert_1^2,\quad \forall t>0.
    \end{equation}
  Then the following Nash type inequality is satisfied
  \begin{equation}\label{nashtilde}
 \vert\vert  f\vert\vert_2^2  \; {\Lambda}(\ln  \vert\vert  f\vert\vert_2^2)
  \leq
  (Af,f),\quad \forall f\in {\mathcal D}(A),\;   \vert\vert 
f\vert\vert_1\leq 1,
  \end{equation}
  where ${\Lambda}(s)=\sup_{t>0}(st-t\ln m(1/2t)\,)$, $s\in \R$.
 \end{Th}
 
  The function ${\Lambda}$ is the conjugate function (or so-called Legendre transform) of
  $t\rightarrow t\ln m(1/2t)$.
   The function $x\rightarrow x{\Lambda}(\ln x)$ is nothing but the function
 ${\tilde {\Theta}}( x)= \sup_{t>0}(x/2t)\log (x/m(t))$ of Proposition II.2 of \cite{c}. The interest of the formulation (\ref{nashtilde}) is that it expresses $ {\tilde {\Theta}}$ in terms of the Legendre transform
 ${\Lambda}$ and is well-appropriate to deal with fractional powers of $A$. Indeed, in \cite{bema}
 it is proved that if $A$ satisfies a Nash type inequality with  ${\tilde {\Theta}}( x)= 
 x{\Lambda}(\ln x)$ then $A^{\alpha}$ ($0<\alpha<1$) satisfies a Nash type inequality  with (roughly)
  ${\tilde {\Theta}}_{\alpha}( x)=x{\Lambda}^{\alpha}(\ln x)$ where 
  ${\Lambda}^{\alpha}=\exp(\alpha\ln {\Lambda})$.
\\
 
\Dem{}  Nash type inequality  (\ref{nashtilde}) is proved by using a 
convexity argument and optimization over the time parameter.  By Jensen's inequality,
$$
e^{-2t(Af,f)/\vert\vert  f\vert\vert_2^2}
\leq
\int_0^{+\infty} e^{-2t{\lambda}}d{\mu}({\lambda})=
\vert\vert  T_tf\vert\vert_2^2/\vert\vert  f\vert\vert_2^2
$$
where  $d{\mu}({\lambda})=d(E_{\lambda}f,f)/\vert\vert  f\vert\vert_2^2$
 and $(E_{\lambda})$ is the spectral decomposition of $A$.
 By assumption (\ref{ultracoul}) and $\vert\vert  f\vert\vert_1\leq 1$, we deduce
 $$
e^{-2t(Af,f)/\vert\vert  f\vert\vert_2^2}
\leq
m(t)/\vert\vert  f\vert\vert_2^2.
$$
This can be written as 
$$
-2t (Af,f)\leq
\vert\vert  f\vert\vert_2^2\left(
\ln m(t)-\ln \vert\vert  f\vert\vert_2^2\right)
$$
or equivalently by changing t by $1/2t$
$$
\vert\vert  f\vert\vert_2^2\left(t\ln \vert\vert  f\vert\vert_2^2
-t\ln m(1/2t)\right)
\leq (Af,f).
$$
We finishes the proof by optimizing over $t>0$. 
\\

 It is easily proved that   Nash type inequality  (\ref{nashtilde}) 
is equivalent  to what we shall call  the relaxed Nash type inequality below
  \begin{equation}\label{nashrelx1}
  \vert\vert  f\vert\vert_2^2\log( \vert\vert  f\vert\vert_2)
  \leq t(Af,f)+\log(\sqrt{m(t)})  \vert\vert  f\vert\vert_2^2,\; 
\forall t>0,  \vert\vert  f\vert\vert_1\leq 1.
   \end{equation}
 
 We compare this inequality with the one we can deduce from 
Davies-Simon Theorem recalled in (\ref{gross}). 
 Under  the assumption of  ultracontractivity property of  
Davies-Simon Theorem and with 
Lemma \ref{convex},  we get the following relaxed Nash type 
inequality :
  \begin{equation}\label{nashrelx2}
 \vert\vert  f\vert\vert_2^2\log( \vert\vert  f\vert\vert_2)
  \leq t(Af,f)+\ M(t)\vert\vert  f\vert\vert_2^2,\; \forall t>0,  
\vert\vert  f\vert\vert_1\leq 1,
 \end{equation}
 and the two functions  $M(t)$ and 
$\log(\sqrt{m(t)})$ from (\ref{ultraM}) and (\ref{ultracoul})  are 
the same i.e
 $M(t)= \log(\sqrt{m(t)}) $ (We assume that (\ref{ultraM}) and 
(\ref{ultracoul}) are equalities).
 So the inequalities (\ref{nashrelx1}) and (\ref{nashrelx2}) are the 
same.
\\

We now state the equivalence of relaxed Nash type inequality and 
LSIWP when $Q(f)=(Af,f)$ is a Dirichlet form (see \cite{fu}).  
We apply truncation method as developped in \cite{bcls}. This result 
is essentially contained in \cite{bima}. We give the sketch of the 
proof for completness.

\begin{Th}\label{nashegallsi}
\begin{enumerate}
\item
Assume that   $Q$ is a Dirichlet form and that the following 
inequality  is satisfied
 \begin{equation}\label{nashrelx3}
 \vert\vert  f\vert\vert_2^2\log( \vert\vert  f\vert\vert_2)
  \leq tQ(f)+\ M(t)\vert\vert  f\vert\vert_2^2,\; \forall t>0,  
\vert\vert  f\vert\vert_1\leq 1.
 \end{equation}
Then 
 \begin{equation}\label{nashrelx4}
 \int  f^2\log_+ (f/ \vert\vert  f\vert\vert_2^2) \,d\mu
   \leq t Q(f)+\tilde{M}(t)\vert\vert  f\vert\vert_2^2,\; \forall t>0
 \end{equation}
 with $\tilde{M}(t)=\frac{1}{c_1}M(t/c_1)+c_2$
 (The constants $c_1,c_2>0$ do not depend on $f$ and  $Q$).
 \item
 Conversely (without any assumptions on the quadratic form $Q$), if (\ref{nashrelx4}) is satisfied then (\ref{nashrelx3}) 
is also  satisfied with $M(t)=\tilde{M}(t)$ for all $t>0$.
 \end{enumerate}
\end{Th}

\Dem{}
\begin{enumerate}
\item
The arguments are taken from \cite{bima} and \cite{bcls}. 
We give the ingredients of the proof.
Let $f$ such that $0\leq f$ and $ \vert\vert  f\vert\vert_2=1$. Let 
$k\in \Z$, we define 
$f_k=(f-2^k)\wedge 2^k$.  Fix $t>0$. By assumption, 
$$
 \vert\vert  f_k\vert\vert_2^2\log( \vert\vert  f_k\vert\vert_2/ 
\vert\vert  f_k\vert\vert_1)
  \leq tQ(f_k)+\ M(t)\vert\vert  f_k\vert\vert_2^2.
  $$
  \vskip 5mm
   We have $2^{k-1}\leq \vert\vert  f_k\vert\vert_2/\vert\vert  
f_k\vert\vert_1$
   and
    $
  2^{2k}\mu(2^k\leq f)
  \leq
  \vert\vert  f_k\vert\vert_2^2
  $.
We set   
$W(g)=tQ(g)+M(t)\vert\vert  g\vert\vert_2^2$ (note that $W$ is also a Dirichlet 
form).
Then we deduce for any $k\in \Z$,
$$
  2^{2k}\mu(2^k\leq f)
  \log 
  2^{k-1}
  \leq 
  W(f_k).
  $$
 By discretisation of the integral,
 $$
 \int f^2\log_+ f\,d\mu
 \leq
 \sum_{k=0}^{\infty} 2^{2(k+1)}\log {2^{k+1}}\mu(2^k\leq f\leq 
2^{k+1}).
 $$
 Hence,
altogether,  we get for some   $c,c'>0$,
$$
   \int f^2\log_+ f\,d\mu
 \leq
c\sum_{k=0}^{\infty}  W(f_{k-1})+
c'
\vert\vert  f\vert\vert_2^2
 $$
 We conclude by the fact that $W$ is a Dirichlet form then
 $$
 \sum_{k=0}^{\infty}  W(f_{k-1})
 \leq
 \sum_{k\in\Z}  W(f_{k})
 \leq
 W(f)
 $$
For a demonstration of the last  statement see \cite{bima}. This finishes the proof of the first statement.

\item
 We apply Lemma \ref{convex}. 
 \end{enumerate}
 
This finishes the proof of this theorem.
\\
 
Combining  Theorem \ref{coulhon}  and  Theorem \ref{nashegallsi}, we obtain  a new proof of LSIWP from ultracontractivity 
property. Indeed, we first  apply Coulhon's result \ref{coulhon} 
then we get the so-called  relaxed Nash type inequality (\ref{nashrelx1}).  We now   apply Theorem
\ref{nashegallsi} to conclude.  Note that in \cite{d}, they get  
$\tilde{M}(t)=M(t)$ for all $t>0$.

\section{Relations between Nash type inequality  and LSIWP}\label{relnashlsi}
\setcounter{equation}{0}

We now prove that Nash type inequality  and LSIWP are (essentially) equivalent. We obtain this result by putting together   arguments of Section \ref{nproflsi} and using  the fact that the function ${\Lambda}$  in   (\ref{nashls})  is a Legendre transform. The natural assumption on ${\Lambda}$ comes from the following remark:
we have proved   that ultracontractivity property  (see (\ref{nashtilde})) or LSIWP 
(see (\ref{lsiwpcoul})  and  (\ref{lsiwpcoul2})) implies a Nash type inequality  of the form
\begin{equation}\label{nashls}
\|f\|^2_2\;{\Lambda}(\ln \|f\|^2_2)\leq (Af,f),  \|f\|_1\leq 1, f\in {\mathcal D}(A).
\end{equation}
where ${\Lambda}$ is given by the  Legendre transformation (i.e  the so-called conjugate function) of the function $\psi$,
$$
 {\Lambda}(y)=\sup_{t>0}(ty/2-\psi(t))
$$
with $\psi(t)=t\beta(1/t)$
(see (\ref{nashtilde})). It implies
 for any $t>$ and any  $y\in \R$,
\begin{equation}\label{leg}
 ty/2-{\Lambda(y)}\leq \psi(t)
\end{equation}
 which says in particular that the function $y\rightarrow  ty/2-{\Lambda(y)}$ is bounded for any $t>0$. Note also, in the theory of Orlicz spaces, the functions  $\psi$ and ${\Lambda}$ are $N$-functions in the sense of \cite{a} and 
 $\psi$ is obtained by the duality  formula, for any $t>0$,
 $$
 \sup_{y\in \Ri}(ty/2-{\Lambda(y)})=\psi(t)
 $$
Here, $\psi$ and ${\Lambda}$ are not necessarily $N$-functions. In our  case, it is not a problem because we only need the inequality (\ref{leg}) in  our applications.
 So as we can see,  Legendre transform (or more generally  convexity) plays again an important role in our theory.
 \\
 
 We are now in a position to give the  natural condition on ${\Lambda}$ to formulate our first statement:
  Nash  type inequality $\Rightarrow$ LSIWP.

 \begin{Th}\label{naslas}
Assume that (\ref{nashls})  holds true for some function ${\Lambda}$ with $(Af,f)$ a Dirichlet form
 and that ${\Lambda}$ satisfies the following hypothesis: 
for any fixed $t>0$, the function $y\in \R\rightarrow ty/2 - {\Lambda}(y)$ is bounded above. Set
$$N(t)=\sup_{y\in\Ri} (ty/2 - {\Lambda}(y)), \;t>0$$
and
$$ \beta(t)={t}N(1/t).$$
Then,  there exists $c_1,c_2>0$ such that, for all $f\in {\mathcal D}(A)$,
\begin{equation}\label{lsiwpnash}
 \int f^2\log_+( f/\vert \vert   f  \vert \vert_2)\,d\mu \leq
 t(Af,f)+{\tilde \beta}(t)\vert \vert   f  \vert \vert_2^2, \quad \forall t>0  
 \end{equation}
 with ${\tilde \beta}(t)=\frac{1}{2c_1}\beta(\frac{2t}{c_1})+c_2$.
\end{Th}

 In fact, the main point of  Theorem \ref{naslas}  is the existence of the function $N$. 
For instance, the second hypothesis of this theorem is satisfied if the following conditions  (A1) and (A2) below holds true.
 
 (A1): For any $s>0$ and any $y_0\in \R$, the function $y\rightarrow sy-{\Lambda}(y)$ is bounded above on the interval $(-\infty,y_0)$.

(A2): $\lim_{y\rightarrow +\infty} \frac{\Lambda(y)}{y}=+\infty$.

For example, the assumption (A1)  is satisfied if $\Lambda$ is non-negative or  if there exists $y_1\in \R$ such that ${\Lambda}(y)=0$  for all $y\leq y_1$ and $\Lambda$  continuous or  $\lim_{y\rightarrow -\infty}  \Lambda(y)=0$ and $\Lambda$  continuous.

 As we have already mentioned,  if ${\Lambda}$ is given by the  Legendre transformation  of some function i.e
 ${\Lambda}(y)=\sup_{t>0}(ty/2-\psi(t))$ (finite at any point $y\in \R$) then ${\Lambda}$ satisfies immediately the hypothesis of our theorem. Also note  that the transformation  
 $$ t\rightarrow  \beta(t)={t}N(1/t).$$is idempotent.  So  $N(t)={t} \beta(1/t)$.

\Dem{}
By Theorem \ref{nashegallsi}, it is enough to prove (\ref{nashrelx3}).  
As a consequence of the definition of $N(t)$, we have for any $t>0$ and any $y\in \R$, 
$$yt/2-{\Lambda}(y)\leq N(t)$$
or equivalently,
$$
yt/2-N(t)\leq {\Lambda}(y).
$$
Let $y=\ln \vert \vert   f  \vert \vert_2^2$ in the inequality just above and multiply it by $\vert \vert   f  \vert \vert_2^2$. Hence, by  our assumption  (\ref{nashls}), we deduce
$$
 t\vert \vert   f  \vert \vert_2^2\ln \vert \vert   f  \vert \vert_2^2
 -N(t)\vert \vert   f  \vert \vert_2^2\leq (Af,f).
 $$
 We set $t=1/s, s>0$. This yields
 $$
  \vert \vert   f  \vert \vert_2^2\ln \vert \vert   f  \vert \vert_2^2
  \leq s(Af,f)+sN(1/s) \vert \vert   f  \vert \vert_2^2
 $$
 and by definition of $\beta$,
 $$
  \vert \vert   f  \vert \vert_2^2\ln \vert \vert   f  \vert \vert_2^2
  \leq s(Af,f)+\beta(s) \vert \vert   f  \vert \vert_2^2.
 $$
 So, (\ref{nashrelx3}) is proved with $M(t)=\frac{1}{2}\beta(2t)$.  Now, we apply Theorem \ref{nashegallsi}.
 This  finishes the proof.
 \\

We now state the converse of Theorem \ref{naslas}. Note that we do not need any assumption on the quadratic form $(Af,f)$ for this converse.

 \begin{Th}\label{naslasconv}
Assume that,  for all $f\in {\mathcal D}(A)$,
\begin{equation}\label{lsiwpnashconv}
 \int f^2\log_+( f/\vert \vert   f  \vert \vert_2)\,d\mu \leq
 t(Af,f)+{ \beta}(t)\vert \vert   f  \vert \vert_2^2, \quad \forall t>0  
 \end{equation}
 is satisfied.
Set
$${\Lambda}(y)=\sup_{t>0} (ty/2 - N(t)), \;y\in \R$$
with
$$N(t)={t} \beta(1/t).$$
Then,  for all $f\in {\mathcal D}(A)$,
\begin{equation}\label{nashlsconv}
\|f\|^2_2\;{\Lambda}(\ln \|f\|^2_2)\leq (Af,f),  \|f\|_1\leq 1.
\end{equation}
\end{Th}

The function ${\Lambda}$ is automatically defined as  can be seen in the course of the proof.
 
 \Dem{}
 The assumption (\ref{lsiwpnashconv}) implies obviously LSIWP. We now repeat the argument of the beginning of the proof of Theorem \ref{anotherbound}. We apply Lemma \ref{convex} to get 
 (\ref{lsiwpcoul4}) with $B=\Lambda$. This completes the proof.
 \\
 
 We have shown that Nash type inequality are   equivalent to LSIWP in the sense of Theorem \ref{naslas} and Theorem \ref{naslasconv}. 

\section{Nash type inequality}\label{nti}
\setcounter{equation}{0}
In this section, we  study  some aspects  of the relationship between 
ultracontractivity and  Nash
type inequality (see \cite{c}, \cite{t}  ) of the form :

\begin{equation}
B(\|f\|^2_2)\leq Q(f)
\end{equation}
for all $f \in Quad(A)\cap L^1$ with $\|f\|_1\leq 1$.
The classical Nash inequality corresponds to $B(t)=ct^{1+2/n}$ (see 
\cite{cks}) i.e.
\begin{equation}\label{nash}
c\|f\|^{2+4/n}_2\leq Q(f)\|f\|^{4/n}_1
\end{equation}

 Let $V,W$ be  two continuous functions on $[0,+\infty[$ and $b_{1}$ 
a continuous
function on $]0,+\infty[$. We begin  by the following easy but 
important proposition :

\begin{pro}\label{nashequiv}
The two following inequalities are equivalent :
for all $t>0$ and all $f \in Quad(A)\cap L^1$ with $\|f\|_1\leq 1$,

\begin{equation}\label{para}
V(\|f\|^2_2)\leq tQ(f)+b_1(t)W(\|f\|^2_2)
\end{equation}

and 

\begin{equation}\label{bin}
B(\|f\|^2_2)\leq Q(f)
\end{equation}
where is defined by $B(x)=\sup_{s>0} (sV(x)-b(s)W(x)), x>0$ with $b(s)=sb_1(1/s)$.
\end{pro}

The    inequality  (\ref{para})  will be called relaxed Nash type inequality. 
We shall apply this proposition in two important cases :
\\

{\underline{case (A)}} :
\begin{equation}
\qquad V(x)=x,\quad \quad W(x)=1.
\end{equation}

{\underline{case (B)}} :
\begin{equation}
\qquad V(x)=\frac{x}{2}\log x,\quad \quad W(x)=x.
\end{equation}

\begin{rem}
The assumptions on $V,W,b_1$ and (\ref{para}) implies that 
$B(\|f\|^2_2)$ is finite.
With $V(x)=x$, $W(x)=1$ then 
$B(x)=\sup_{s>0}(sx-b(s))$ is the complementary function
(or Legendre transform) of $b$.
In applications, we often recover $b(s)$ by the same formula 
$b(s)=\sup_{x>0}(xs-B(x))$ (see \cite{a} p.229). 

With the choice $V(x)=\frac{x}{2}\ln x$ and $W(x)=x$,
$B(x)=x\sup_{s>0}(\frac{s}{2}\ln x-b(s))$ is similar to the function ${\tilde \Theta}$ introduced in section \ref{nproflsi}. 
\end{rem}

\Dem {}
 The proof is easy and left to the reader.
 \\

 By the remark just above,
 the problem of  obtaining  the Nash type inequality (\ref{bin}) is 
equivalent to show the inequality
(\ref{para}). But some difference may arise  from the choice of the functions $V$ and $W$ as we shall see below.
Under an ultracontractivity assumption on the semigroup, we apply 
the proposition
\ref{nashequiv} in case (A) and (B) respectively.
\\

{\underline{Case A}}

We start with this simple case.

\begin{Th}\label{nacla}
Let $T_t$ be an ultracontractive semigroup such that
\begin{equation}\label{reg}
\|T_tf\|_{\infty}\leq a(t)\|f\|_1
\end{equation}
then 
\begin{equation}\label{superpoin}
\|f\|_2^2\leq tQ(f)+a(t)
\end{equation}
for all $t>0$ and $f\in Quad(A)\cap L^1$, with $\|f\|_1\leq 1$.
\\

and
\begin{equation}\label{nashsuperpoin}
B(\|f\|^2_2)\leq Q(f)
\end{equation}
for all $t>0$ and $f\in Quad(A)\cap L^1$, with $\|f\|_1\leq 1$,
where we set  $b(s)=sa(1/s)$ and $B(x)=\sup_{s>0} (sx-b(s))$.
\end{Th}

The inequality  (\ref{superpoin}) is called super-Poincar\'e  inequality in \cite{w} (see (1.2) of   \cite{w}). See also Sec.5 of \cite{w}  for related results to the ultracontractivity property. 

\Dem {}
Let $f\in L^1\cap Quad(A)$ with $\|f\|_1\leq 1$. Set $f_s=T_sf$, then 
for all $t>0$ :
$$ f=f_t+\int_0^t Af_s\,ds $$
Hence,
$$\|f\|_2^2=(f_t,f)+ \int_0^t (Af_s,f)\,ds$$
and
$$\|f\|_2^2\leq \|f_t\|_{\infty}\|f\|_1+\int_0^t Q(f_{s/2})\,ds$$
The function $s\mapsto Q(f_s)=(Af_s,f_s)$ is non-increasing, thus
$$\|f\|_2^2\leq a(t)+tQ(f)$$
We conclude the proof by applying Proposition \ref{nashequiv} 
with $V(x)=x, W(x)=1$ and 
$b_1(t)=a(t)$.
\\

In particular, if $a(t)=ct^{-\frac{n}{2}}$ then we obtain the 
classical Nash inequality  (\ref{nash}).
This result is well-known (see \cite{vsc}). In fact,   (\ref{reg}) and
(\ref{nash}) are equivalent (see\cite{cks}) and are also equivalent 
to Sobolev
inequality (when $n>2$) :

\begin{equation}
\|f\|^2_{\frac{2n}{n-2}}\leq cQ(f)
\end{equation}
This last result is due to Varopoulos (see 
\cite{v}).
\\

For  applications,  an  important family of functions is
\begin{equation}\label{at}
 a(t)=c_1t^{-\alpha}\exp(\frac{c_2}{t^{\gamma}}),\quad{\alpha}\geq 0, {\gamma}\in {\R}.
 \end{equation}
 This kind of function motivates our study. Of course, when $t$ is small the main term is the exponential for which the computation of $B$ just above is rather complicated.
Indeed, with  ${\gamma}\not=0$ and ${\alpha}\not=0$ the function  $B$ doesn't 
seem to be known explicitely. We can certainly estimate $B(x)$ when $x$ is large.
 But we shall see below that a better approach of Nash type inequality is to
 use the logarithmic Sobolev inequality. In that case $a(t)$ is replaced by $\ln a(t)$  with a change of 
 the couple $(V,W)$ and the corresponding function $B(t)$ can be computed exactly (see case B).
 \\

{\underline{Case B}} :

We have at least  two possibilities  to prove  Nash type inequality  from an ultracontractivity property 
for the case   $V(x)=\frac{x}{2}\ln x$ and $W(x)=x$. One way is to use the LSIWP and  
 the convexity Lemma \ref{convex}. An alternative proof is to use Coulhon's result recalled in Theorem \ref{coulhon}.

\begin{Th}\label{nashlog}
Let $T_t$ be an ultracontractive semigroup such that
\begin{equation}\label{reg1}
\|T_tf\|_{\infty}\leq a(t)\|f\|_1
\end{equation}
or
\begin{equation}\label{reg2}
\|T_tf\|_{\infty}\leq \sqrt{a(t)}\,\|f\|_2.
\end{equation}
Then, for all $t>0$ and $f\in Quad(A)\cap L^1$, with $\|f\|_1\leq 1$ :

\begin{equation}\label{iks}
\|f\|^2_2\log \|f\|_2 \leq tQ(f)+\frac{1}{2}\log a(t)\|f\|^2_2
\end{equation}
Or equivalently the following Nash type inequality 
\begin{equation}\label{lognor}
B(\|f\|^2_2)\leq Q(f)
\end{equation}
with $B(x)=x\sup_{s>0}(\frac{s}{2}\log x -b(s))$ and
 $b(s)=\frac{s}{2}\log a(1/s)$.
\end{Th}

\Dem {}
By applying Theorem \ref{gross}  with $M(t)=\frac{1}{2}\log 
a(t)$,  the  following  LSIWP  is satisfied: for all $t>0$, 
\begin{equation}
\int f^2\log (f/ \|f\|_2) \, d{\mu} \leq
t Q(f) +M(t) \|f\|^2_2 
\end{equation}

From the lemma \ref{convex}, we deduce (\ref{iks}). 
We conclude by applying  Proposition \ref{nashequiv} with 
$V(x)=\frac{x}{2}\log x, W(x)=x$ 
and $b_1(t)=\frac{1}{2}\log a(t)$. This completes the proof.
\\

The  alternative proof using  Coulhon's result is left to the reader
(see Theorem \ref{coulhon}).
\\

We now give an example.
When ${\alpha}=0$ in (\ref{at}), we are able to explicit  function $B$ in  
Nash type inequality.

\begin{Th}\label{corexp}
If the following inequality is satisfied for some  ${\gamma}>0$,
\begin{equation}\label{e1}
\|T_tf\|_{\infty}\leq c_1\exp(\frac{c_2}{t^{\gamma}})\|f\|_1
\end{equation}
Then there exist $k,{\beta}>0$ such that for all $f\in Quad(A)\cap 
L^1$, with $\|f\|_1\leq 1$,
\begin{equation}\label{e2}
k\|f\|^2_2\left[\left(\log 
\frac{\|f\|_2}{\beta}\right)_+\right]^{1+\frac{1}{\gamma}} \leq
Q(f)
\end{equation}
Conversely (\ref{e2}) implies (\ref{e1}) with different constants 
$c_1$ and $c_2$.
\end{Th}

\Dem {}
We apply Theorem \ref{nashlog}.
Since $B(t)=tD(\log \sqrt{t})$, we have just  to compute 
$D(x)=\sup_{s>0}(sx-b(s))$ with $b(s)=\frac{\log 
c_1}{2}s+\frac{c_2}{2}s^{1+{\gamma}}$. 
For $x\in \R$, we easily study the extrema of 
$h(s)=sx-b(s)$ and obtain
$D(x)=k\left[(x-k_1)_+\right]^{1+\frac{1}{\gamma}}$, $x\in \R$ with 
$k>0$. Let ${\beta}=e^{k_1}$
then  $B(t)=kt\left[\left(\log 
\frac{\sqrt{t}}{\beta}\right)_+\right]^{1+\frac{1}{\gamma}}$.
The converse is proved by applying  Proposition \ref{nashequiv}   
and Theorem \ref{coroda}.
The proof is completed.
 
\begin{rem}
We denote that $b(s)$ may not an N-function in the sense of \cite{a} 
p.228 since $\lim_{s\rightarrow 0} b(s)/s$ is not necessarily zero.
But for all $t,s>0$, the functions $B$ and $b$ satisfy
$$
t(s\log \sqrt{t} - b(s))\leq B(t)
$$
Let $v=\log \sqrt{t}$, then
$$
sv-H(v)\leq b(s)
$$
with
$H(v)=e^{-2v}B(e^{2v})$. Hence,
$$
{\tilde H}(s)=\sup_{v\in \R}(sv-H(v))\leq b(s)
$$
It is easily  seen that ${\tilde H}(s)=b(s)$ when 
$B(t)=kt\left[\left(\log 
\frac{\sqrt{t}}{\beta}\right)_+\right]^{1+\frac{1}{\gamma}}$.
\end{rem}

We now show that an ultracontractive bound on the semigroup implies 
a generalized Gross' inequality (see (\ref{ined1} below). Such an equality  is  easily deduced from LSIWP. 
 
\begin{Th}\label{betnash}
Let $T_t$ be an ultracontractive semigroup such that
\begin{equation}
\|T_tf\|_{\infty}\leq a(t)\|f\|_1
\end{equation}
then
\begin{equation}\label{ined1}
D\left(\int f^2\log f\,d{\mu}\right) \leq
Q(f)
\end{equation}
for all $0\leq f\in Quad(A)\cap L^1\cap L^{\infty}$ with  $\|f\|_2=1$ 
and where 
$D(t)=\sup_{s>0}(st-b(s))$ with $b(s)=\frac{s}{2}\log a(\frac{1}{s})$.
Equivalently,
\begin{equation}\label{ined2}
\|f\|_2^{_2} \; D\left(\|f\|_2^{-2} 
\int f^2\log \left(f/\|f\|_2\right) \,d{\mu}\right) \leq
Q(f)
\end{equation}
for all $0\leq f\in Quad(A)\cap L^1\cap L^{\infty}$.
\end{Th}

\Dem {}
By Davies-Simon Theorem   \ref{gross}, we have 
\begin{equation}
\int f^2\log f \, d{\mu}\leq
t Q(f) +\frac{1}{2}\log a(t) \|f\|^2_2 + \|f\|^2_2\log \|f\|_2
\end{equation}
 for all $t>0$.
We now assume that $\|f\|^2=1$ then changing $t$ by $1/s$,

\begin{equation}
s\left(\int f^2\log f\, d{\mu}\right) -\frac{s}{2}\log a(1/s) 
\leq Q(f) 
\end{equation}

Thus, with $D$ define as above, we conclude (\ref{ined1}) and 
(\ref{ined2}) by renormalisation.
This proves the theorem.
\\

We can see Theorem  \ref{nashlog} as a corollary of Theorem 
\ref{betnash} when $D(t)$ is
non decreasing in $t$. The proof is a simple application of the lemma 
\ref{convex} which can be
also formulated as follows, for all $0\leq f\in L^1\cap L^{\infty}$ with 
$\|f\|_1\leq 1$ :
\begin{equation}
\log \|f\|_2
\leq 
\|f\|_2^{-2} 
\int f^2\log \left(f/\|f\|_2\right) \,d{\mu}
\end{equation}
then from (\ref{ined2}),
\begin{equation}
\|f\|_2^2 D(\log \|f\|_2)
\leq 
Q(f)
\end{equation}
But the first term of this inequality is $B(\|f\|_2^2)$ in 
(\ref{lognor} ) because
$B(t)=tD(\log \sqrt{t})$. We conclude Theorem  \ref{nashlog}.
\\

\section{Weak Sobolev inequalities and Nash Type inequalities. }
\setcounter{equation}{0}

We apply  the results of the preceding  section with  
$a(t)=c_1t^{-\frac{n}{2}}$  and  recall   explicitely some well-known results 
concerning semigroups with  this polynomial ultracontractivity.  Recall that, in our setting,
\begin{equation}
\| e^{-At}\|_{\infty,2}\leq c_1t^{-\frac{n}{4}}
\end{equation}
is equivalent  to 
\begin{equation}
\| e^{-At}\|_{\infty,1}\leq c'_1t^{-\frac{n}{2}}
\end{equation}
Let ${\cal D}$ be the domain of $A$.

\begin{Th}
The following inequalities are equivalent

(i)
\begin{equation}\label{reg1bis}
\| e^{-At} f\|_{\infty}\leq c_1t^{-\frac{n}{4}} \vert\vert f \vert\vert_2
\end{equation}
for all $f\in L^2$ and for all $t>0$.

(ii)
\begin{equation}\label{ws1}
\int f^2\log f \, d{\mu}\leq
\log \left( c_3 Q(f)^{\frac{n}{4}}\right)
\end{equation}
for all $f\in {\cal D}\cap L^1\cap L^{\infty}, f\geq 0$ and 
$\|f\|_2=1$.

(iii)
\begin{equation}\label{n1}
\|f\|_2^{2+\frac{4}{n}}\leq c_4 Q(f)\|f\|_1^{\frac{4}{n}}
\end{equation}
for all $f\in L^1\cap {\cal D}$.

(iv)
\begin{equation}\label{s1}
\|f\|_{\frac{2n}{n-2}}^2\leq
c_5 Q(f)
\end{equation}
for all  $f\in {\cal D}$ (when  $n>2$).

\end{Th}

\Dem {}
The equivalences  between (\ref{reg1bis}) and (\ref{n1}) and (\ref{s1}) 
(when $n>2$) are  well-known, see
\cite{vsc},\cite{d}.
\\

In \cite{b}, D.Bakry gives the arguments to prove
that (\ref{ws1}) is equivalent to (\ref{reg1bis})  (see  remark  at 
the end of
page 64  and Section 5  p.67 of \cite{b})).  
Here, we focus on the inequality (\ref{ws1}).
We first prove that (\ref{reg1bis}) implies (\ref{ws1}) by applying
Theorem \ref{betnash} with
$a(t)=c_0\,t^{-\frac{n}{2}}$ .  A simple computation gives us
$D(y)=ce^{\frac{4y}{n}},\, y\in \R$. 
Then, we deduce  
for all $f\in {\cal D}\cap L^1\cap L^{\infty}, f\geq 0$ and 
$\|f\|_2=1$ :

\begin{equation}
c\,\exp\left(\frac{4}{n}\int f^2\log f\, d{\mu}\right) \leq Q(f)
\end{equation}
because  $D$ is increasing  and  thus invertible, this inequality
 is equivalent to (\ref{ws1}) i.e
\begin{equation}
\int f^2\log f\, d{\mu} \leq
\log \left( c_3 Q(f)^{\frac{n}{4}}\right)
\end{equation}
with 
$c_3=c^{-\frac{n}{4}}$.
\\

We now prove that (\ref{ws1}) implies (\ref{n1})  by applying
Theorem  \ref{nashlog}. 
 We get, for any  $0\leq f\in {\cal 
D}\cap L^1\cap L^{\infty}$,
\begin{equation}
\|f\|_2^{1+\frac{n}{2}}\leq c_3 Q(f)^{\frac{n}{4}}
\end{equation}
with  $\|f\|_1\leq 1$.  This proves (\ref{n1}).
\\

The original proof of the implication   
(\ref{reg1bis})  to  (\ref{n1}) is obtained by using Theorem  (\ref{nacla}).
But, as we have seen above, the applications of  Theorem \ref{nacla}
have  some limitations. So we shall prefer the approach with LSIWP given  by Theorem \ref{betnash}.
\\

The inequality (\ref{ws1}) is called weak-Sobolev inequality of dimension
$n$  in  \cite{b}.  In section \ref{nti}, this  inequality is generalized by inequality (\ref{ined1}).
Then (\ref{ws1}) is a particular case of  (\ref{ined1}) with  the function 
 $D(y)=ce^{4y/n},\, y\in \R$ which is invertible. More generally when $D$ is invertible,
 the inequality  (\ref{ined1}) is equivalent to 
\begin{equation}\label{wsd}
\int f^2\log f \, d{\mu}\leq
D^{-1}\left( Q(f) \right)
\end{equation}
for all $f\in {\cal D}\cap L^1\cap L^{\infty}, f\geq 0$ and 
$\|f\|_2=1$.

\section{Examples on the infinite Torus}\label{sectorus}
\setcounter{equation}{0}
 In this section, we give  families of examples of semigroups with one-exponential and double-exponential ultracontractivity property when the time is small (see definitions below).
  In fact, a 
natural
setting for having such behaviors is the infinite torus  $\T^{\infty}$. The reference for such examples is the paper by A.Bendikov \cite{b2} (see also \cite{b1}). This paper also contains  much more  examples of classes of ultracontractivity. In fact, much  weaker behavior than polynomial ultracontractivity  can be produced (see\cite{b2}).
\\

Let  $\T^{\infty}$ be the infinite  torus with its ordinary product 
structure. The group 
$\T^{\infty}$ is a compact abelian group. The neutral element is 
denoted by
$0$ and  $dm$ the normalised Haar measure. This measure  is the
countable product of the normalised Haar measure of $\T$ which is 
identified
with  $[-{\pi},{\pi}]$.
Let ${\mu}_t$ the brownian semigroup on $\T$ and ${\cal
A}=\{a_k\}_{k=1}^{\infty}$ a sequence of strictly positive numbers.
For each $t>0$, we define  the product measure:
$${\mu}_{t}^{\cal A}= \otimes_{k=1}^{\infty} {\mu}_{ta_k}$$
then $({\mu}_{t}^{\cal A})_{t>0}$ defines  a symmetric convolution
semigroup on  $\T^{\infty}$ denoted by $(T_t^{\cal A})_{t\geq0}$  and
$$|| T^{\cal A}_t||_{1\rightarrow {\infty}}={\mu}_{t}^{\cal A}(0).$$
One  important aspect  of such semigroups for  Harmonic Analysis  theory in infinite
dimensional spaces is that  $({\mu}_{t}^{\cal A})$ is not necessarily
absolutely continuous with respect to the Haar measure $dm$.
We have to impose conditions on the sequence ${\cal A}$ in order to 
get a
continuous density.
To this purpose,  we set  for $x>0$ : 
$${\cal N}^{\cal A}(x)=\sharp\{ k\geq 1: a_k\leq x\}.$$

 It is proved in \cite{b2} (Thm 3.6 p.51)  that if 
$$\log {\cal N}^{\cal A}(x)=o(x)\quad \hbox {as }\quad x\rightarrow 
+\infty$$
 then ${\mu}_{t}^{\cal A}$ has a continuous density. The converse is also true.  
We also denote  by  ${\mu}_{t}^{\cal A}$ the density when it exists. The infinitesimal generator is given formally as an infinite Laplacian $A=\sum_{k=1}^{\infty} a_k\frac{\partial^2}{\partial^2 x_k}$.
In the two following subsections, we focus on two particular examples of  ultracontractivity property.

\subsection{One-exponential ultracontractivity }
We shall say that a semigroup $(T_t)$  satisfies a one-exponential ultracontractivity  property at zero if there exists $\alpha>0$ and $t_0>0$ such that
 \begin{eqnarray}\label{expobehaup}
||T_t||_{1\rightarrow +\infty}
\leq
c'\exp\left(\frac{c}{t^{\alpha}}\right), \quad 0<t<t_0.
 \end{eqnarray}
for some constants $c,c'>0$.
We shall say that a semigroup $(T_t)$ has a strict one-exponential ultracontractivity property if 
(\ref{expobehaup}) is satisfied and moreover
 \begin{eqnarray}\label{expobehalow}
c'\exp\left(\frac{c}{t^{\alpha}}\right)\leq
||T_t||_{1\rightarrow +\infty},
 \quad 0<t<t_0.
 \end{eqnarray}
 with the same index $\alpha$ but with possibly other constants $c,c'$.
\\

Note that  of (\ref{expobehaup})  holds for some $t<t_0$ then it holds for any $t>0$.
Indeed, the function $||T_t||_{1\rightarrow +\infty},$ is non-increasing in $t$.
We have the following family of examples satisfying a  strict one-exponential ultracontractivity
property:
\begin{Th}\label{oneben} ( \cite{b2} Thm 3.18 p. 54).
Let ${\alpha} >0$.
Assume that  ${\cal N}^{\cal A}(x)\sim x^{\alpha}$ as $x\rightarrow +\infty$.
Then we have
 \begin{eqnarray}
\log {\mu}_t^{\cal A}(0) \sim k(\alpha)\; t^{-{\alpha}}\hbox{ as}\quad
t\rightarrow  0
 \end{eqnarray}

In particular, there exists $t_0$,  $c_1(\alpha ), c_2(\alpha )$ such 
that :
$\forall t\in ]0,t_0]$ ,
 \begin{eqnarray}\label{onexp}
\exp\left(\frac{c_1({\alpha})}{t^{\alpha}}\right)
\leq 
||T_t||_{1\rightarrow +\infty}
\leq
\exp\left(\frac{c_2({\alpha})}{t^{\alpha}}\right)
 \end{eqnarray}
\end{Th}

Recall that
 $||T_t||_{1\rightarrow +\infty}={\mu}_t^{\cal A}(0)$. The condition on the sequence
 ${\cal A}$ is satisfied,
 for instance with  $a_k=k^{1/{\alpha}}, k\geq 1$ with $\alpha>0$. 
 \\
 
Now, we  apply Theorem \ref{corexp}  
 with $X=\T^{\infty}$
and ${\mu}$ the Haar measure on $X$.
We have the following Nash type inequality.

\begin{Th} 
Under the assumptions of Theorem  \ref{oneben}. There exists  constants  $c_1,c_2>0$ such that :
\begin{eqnarray}\label{normlogalpha}
||f||_2^2\left[\log_+ \left(\frac{|| f||_2}{c_1|| f||_1}\right)\right]^{1+{1\over
{\alpha}}}
\leq  c_2 Q(f)
\end{eqnarray}
for all $0\leq f\in  {\cal D}\cap L^1\cap
L^{\infty}$.
\end{Th}

We deduce the folllowing result by using cut-off method developed  as in \cite{bcls}  (see  \cite{bima} and also \cite{w} for such inequalities).

\begin{Th}
There exists $c_3,c_4>0$ such that,
for all $0\leq f\in  {\cal D}\cap L^1\cap \in
L^{\infty}$,
 \begin{eqnarray}\label{logalpha}
\int_{\T^{\infty}} f^2 \left[\log\left(\frac{f}{c_3|| 
f||_2}+1\right)\right]^{1+{1\over
{\alpha}}}d\,{\mu}
\leq
c_4( Q(f)+|| f|_2^2).
 \end{eqnarray}
\end{Th}

Conversely, by a convexity argument  similar to Lemma \ref{convex}  (in the genral framework), (\ref{logalpha}) implies 
(\ref{normlogalpha}) with $Q(f)+|| f|_2^2$ instead of  $Q(f)$.

\subsection{Double-exponential ultracontractivity }\label{ddsect}
We shall say that a semigroup $(T_t)$ satisfies  a double-exponential ultracontractivity  property at zero if there exists $\gamma>0$ and $t_0>0$ such that
 \begin{eqnarray}\label{expexpobehaup}
||T_t||_{1\rightarrow +\infty}
\leq
c'\exp(\exp\left(\frac{c}{t^{\gamma}}\right)), \quad 0<t<t_0,
 \end{eqnarray}
for some constants $c,c'>0$.
We shall say that a semigroup $(T_t)$ has a strict double-exponential ultracontractivity property  if 
(\ref{expexpobehaup}) is satisfied and moreover
 \begin{eqnarray}\label{expexpobehalow}
c'\exp(\exp\left(\frac{c}{t^{\gamma}}\right))\leq
||T_t||_{1\rightarrow +\infty},
 \quad 0<t<t_0,
 \end{eqnarray}
 with the same index $\gamma$ but with possibly other constants $c,c'$.
\\

Note that  if (\ref{expexpobehaup}) holds for some $t<t_0$ then it holds for any $t>0$.
We have the following family of examples satisfying a strict double-exponential ultracontractivity
property:
 
\begin{Th}\label{doubleben} ( \cite{b2} Thm 3.27 p. 59).
Let ${\gamma} >0$.
Assume that $\;\log {\cal N}^{\cal A}(x)\sim x^{\gamma/\gamma+1}$ as $x\rightarrow +\infty$.
Then we have
 \begin{eqnarray}
\log\log {\mu}_t^{\cal A}(0) \sim c(\gamma)\; t^{-{\gamma}}\hbox{ as}\quad
t\rightarrow  0
 \end{eqnarray}

In particular, there exists $t_0,c_1(\gamma ),c_2(\gamma )>0$ such 
that :
$\forall t\in ]0,t_0]$ ,
 \begin{eqnarray}\label{dexp}
\exp(\exp\left(\frac{c_1({\gamma})}{t^{\gamma}}\right))
\leq 
||T_t||_{1\rightarrow +\infty}
\leq
\exp(\exp\left(\frac{c_2({\gamma})}{t^{\gamma}}\right))
 \end{eqnarray}
\end{Th}

The condition on the sequence
 ${\cal A}$ is satisfied,
 for instance  with $a_k=\left[\ln (k+2)\right]^{\delta}, k\geq 1$ with $\delta=
 (\gamma+1)/\gamma$. 
 \\
 
Now, we  apply Theorem \ref{corexp}  
 with $X=\T^{\infty}$
and ${\mu}$ the Haar measure on $X$.
We have the following Nash type inequality,

\begin{Th} 
Under the assumptions of Theorem  \ref{doubleben}. There exists  constants  $c_1,c_2>0$ such that :
\begin{eqnarray}
||f||_2^2 \log\left(\frac{|| f||_2}{|| f||_1}\right)
\left[\log \left( c_1 \log_+ \frac{|| f||_2}{|| f||_1}\right) \right]_+^{{1\over
{\gamma}}}
\leq  c_2Q(f)
\end{eqnarray}
for all $0\leq f\in  {\cal D}\cap L^1\cap 
L^{\infty}$.
\end{Th}

We deduce the folllowing result. We set  $D_{\gamma}(y)=y_+\left[\log_+y_+\right]^{1/\gamma}$

\begin{Th}
There exists $c_3,c_4>0$ such that,
for all $0\leq f\in  {\cal D}\cap L^1\cap \in
L^{\infty}$,
 \begin{eqnarray}
\int_{\T^{\infty}} f^2 
D_{\gamma}\left( c_3 \log\left(\frac{f}{8|| f||_2}\right)\right)
d\,{\mu}
\leq
c_4 Q(f)
 \end{eqnarray}
\end{Th}

\section {The double-exponential 
case}\label{nop}
\setcounter{equation}{0}
In this section, the assumptions are the same as in Section \ref{relultlsi}.
In this general framework, we  consider the relationship  between the double-exponential  ultracontractivity property and LSIWP. We recall that 
the semigroup satisfies the  double-exponential  
ultracontractivity
property if:
\begin{eqnarray}\label{double}
||T_tf||_{\infty}\leq a(t) ||f||_1, \quad t>0
\end{eqnarray}
with  $a(t)=\exp (\exp(/{t^{\gamma}}))$ with 
${\gamma}>0$ and $c>0$ (see Subsection \ref{ddsect}). 
Applying
Theorem \ref{gross}, we get  the following logarithmic Sobolev inequality with parameter : for all $t>0$  ,
\begin{eqnarray}\label{logdouble}
\int f^2\log f\,d{\mu} \leq tQ(f)+{\beta}(t) ||f||_2^2 + 
||f||_2^2\log ||f||_2
\end{eqnarray}
where ${\beta}(t) =\frac{1}{2} \log a(t) =\frac{1}{2}   \exp(c/{t^{\gamma}})$.

We now discuss about the converse result.
Assume that (\ref{logdouble}) holds true.
We remark that  by formula (\ref{integral}) of  Theorem \ref{coroda} gives  $M_{\eta}(t)=\infty$ for all 
$t>0$ and all ${\eta}>-1$.
So we get no information about ultracontractivity property of the semigroup with this formula.
It is not completely surprising. Indeed, there
exists a semigroup satisfying  (\ref{logdouble}) with  
${\gamma=1}$ but which is not 
ultracontractive 
 hence (\ref{ulb}) doesn't not hold
(\cite{d}: example 2.3.5 p.73 and \cite{ds}: (3) p.355 and p.359). 
So LSIWP  (\ref{logdouble}) may not imply the ultracontractivity property of the semigroup 
(at least in that case  when  $\gamma\geq 1$). In that section, we show that, if 
the
condition (\ref{logdouble}) is satisfied with  $\gamma <1$, then the semigroup is ultracontractive of double-exponential type but with $\gamma'$  different  from $\gamma$.
Moreover $\gamma'$  tends to infinity as $\gamma$ tends to 1.

This remark implies  that Nash type inequality and  ultracontractivity property
are not equivalent in general (see \cite{c} and also \cite{bema}).

By the converse Theorem of Davies-Simon \ref{gammap},  if we suppose 
that the
condition (\ref{logdouble}) is fulfilled with ${\gamma}\in ]0,1[$, 
then the semigroup is ultracontractive. But Theorem \ref{gammap} 
doesn't seem to give easily
and explicitely a function $b(t)$ such that : 
\begin{eqnarray}\label{ulb}
||T_tf||_{\infty}\leq b(t) ||f||_1
\end{eqnarray}

In this paragraph, we modify the argument of proof of  Lemma 
\ref{lemphi} to deal with the
double-exponential case. 
\begin{pro}\label{expo}
If (\ref{logdouble})  is satisfied with 
${\beta}(t)=c_1\exp(\frac{c_2}{t^{\gamma}})$,
$0<{\gamma}<1$ then
 (\ref{ulb}) holds with $b(t)=k_1\exp\left(\exp(\frac{k_2}{t^{\gamma'}} )\right)$, 
${\gamma'}=\frac{\gamma}
{1-{\gamma}}$ where $k_1,k_2$ are some positive constants.
\end{pro}

\begin{rem}
\begin{enumerate}
\item
If we apply Theorem \ref{gross} and its (partial) converse 
proposition 
\ref{expo}, we stay in the same  class of functions of type double-exponential. But we lose the exponent  ${\gamma}$.
The question to know if the expression of the exponent ${\gamma'}=\frac{\gamma}
{1-{\gamma}}$ above is optimal is open.
\item
 We also note that ${\gamma'}$  is singular when ${\gamma}$ tends 
to $1$.  By a preceding remark  
 ${\gamma}=1$ is really  a critical index.
\end{enumerate}
 \end{rem}
 
\Dem{}
The proof follows the same lines as the proof of  Theorem \ref{our}
The only change  we need is a modification  of the lemma \ref{lemphi}
when the function $b(t)=c_1\exp(\frac{c_2}{t^{\gamma}} )$.
This is done  with  the following lemma :

\begin{lem}\label{lemphideux}
Suppose that  ${\Phi}(s)\in
C^1([0,+\infty[)$  satisfies ,
  
\begin{equation}\label{fi}
{\Phi}(s)\leq
(-t/2){\Phi}^{\prime}(s) +c_1\exp(\frac{c_2}{t^{\gamma}} )
\end{equation}
for all $s,t>0$ and for a fixed ${\gamma}$ such that $0<{\gamma}<1$. 
Then we have 
\begin{equation}
{\Phi}(t)\leq
k_1\exp(\frac{k_2}{t^{\alpha}} )
\end{equation}
for all $t>0$, where ${\alpha}=\frac{\gamma}
{1-{\gamma}}$, $k_1 =2c_1$ and $k_2= c_2^{\frac{1}{1-{\gamma}}}
\left( \frac{\gamma}{1-{\gamma}}\right)^{\frac{\gamma}{1-{\gamma}}}$.
\end{lem}

{\bf Proof of the lemma :}
Set $t=\frac{ s^{\beta+1} } {\lambda}$ with ${\lambda}>0$ and 
${\beta}>0$ choosen later. We multiply  (\ref{fi}) by $\exp(\frac{-2c_3}{s^{\beta}} ) 
s^{-{\beta}-1}$ with
$c_3=c_2{\lambda}^{\gamma}$. Then ,

\begin{eqnarray}\label{mul}
s^{-{\beta}-1}\exp(\frac{-2c_3}{s^{\beta}} ) 
{\Phi}(s)\leq
\frac{-1}{2{\lambda}}\exp(\frac{-2c_3}{s^{\beta}} ) 
{\Phi}^{\prime}(s) 
\end{eqnarray}

$$
+
c_1s^{-{\beta}-1} 
\exp(\frac{c_3} {s^{\gamma ({\beta}+1)}} -\frac{2c_3}{s^{\beta}} ).
$$
We choose ${\beta}$ such that ${\gamma}({\beta}+1)={\beta}$ then 
${\beta}=\frac{\gamma}{1-{\gamma}}>0$.
We integrate (\ref{mul}) over the interval $[0,t]$ for $t>0$. Let
$I(t)=\int_0^t s^{-{\beta}-1}\exp(\frac{-2c_3}{s^{\beta}}) 
{\Phi}(s)\, ds$. Thus 
\begin{eqnarray}
I(t)\leq  \frac{-1}{2{\lambda}}A(t)+B(t)
\end{eqnarray}
with 
\begin{eqnarray}
A(t)= \int_0^t \exp(\frac{-2c_3}{s^{\beta}} ) 
{\Phi}^{\prime}(s)\, ds
\end{eqnarray}
and 
\begin{eqnarray}
B(t)=c_1\int_0^t s^{-{\beta}-1}\exp(\frac{-c_3}{s^{\beta}} ) 
\, ds.
\end{eqnarray}
These integrals converge because  ${\beta}>0$ and $c_3>0$. The function 
$B(t)$ 
can be explicitely
computed 
\begin{eqnarray}
B(t)=\frac{c_1}{c_3{\beta}}
\exp(\frac{-c_3}{t^{\beta}} ).
\end{eqnarray}
$A(t)$ is computed by integration by parts 
\begin{eqnarray}
A(t)=
\exp(\frac{-2c_3} {t^{\beta}}  ) {\Phi}(t)-2{\beta}c_3 I(t).
\end{eqnarray}
From  (\ref{mul})  we get

\begin{eqnarray}
I(t)\leq \frac{-1}{2{\lambda}}
\exp(\frac{-2c_3}{t^{\beta}}) {\Phi}(t)+\frac{{\beta}c_3}{\lambda}I(t)
+ \frac{c_1}{c_3{\beta}}
\exp(\frac{-c_3}{t^{\beta}}) 
\end{eqnarray}

Now, we chose ${\lambda}$ such that ${\lambda}={\beta}c_3$ then 
${\lambda}=({\beta}c_2)^{\frac{1}{1-{\gamma}} }$. Finally, 
\begin{eqnarray}
\frac{1}{2{\lambda}}\exp(\frac{-2c_3}{t^{\beta}} ) {\Phi}(t)
\leq 
\frac{c_1}{c_3{\beta}}
\exp(\frac{-c_3}{t^{\beta}} ) 
\end{eqnarray}
We conclude the lemma with 
$k_1=2c_1$ and $k_2=c_2^{\frac{1}{1-{\gamma}} }
(\frac{\gamma}{1-{\gamma}} )^{\frac{\gamma}{1-{\gamma}} }$
and then Proposition \ref{expo} is proved.


\end{document}